\documentclass[12pt]{article}
\usepackage{bbm}
\usepackage{amsfonts}
\usepackage{pifont}
\usepackage{hyperref}
\usepackage{mathrsfs}
\usepackage{indentfirst}
\usepackage{amsmath}
\usepackage{amssymb}
\usepackage[amsmath, thmmarks]{ntheorem}
\usepackage{cite}
\usepackage{graphicx}
\usepackage{tikz}
\newtheorem{definition}{\bf Definition}[section]
\newtheorem{lemma}{\bf Lemma}[section]
\newtheorem{theorem}{\bf Theorem}[section]

\newtheorem{corollary}{\bf Corollary}[section]
\newtheorem{proposition}{\bf Proposition}[section]
\newtheorem{example}{\bf Example}[section]

\def\QEDopen{{\setlength{\fboxsep}{0pt}\setlength{\fboxrule}{0.2pt}\fbox{\rule[0pt]{0pt}{1.3ex}\rule[0pt]{1.3ex}{0pt}}}} 
\def\QED{\QEDopen}
\def\proof{{\bf Proof.} }
\def\endproof{\hspace*{\fill}~\QED\par\endtrivlist\unskip}

\begin{document}
\setcounter{page}{1}

\title{{\textbf{On the uniqueness of $L$-fuzzy sets in the representation of families of sets}}\thanks {Supported by National Natural Science
Foundation of China (No.61573240)}}
\author{Peng He\footnote{\emph{E-mail address}: 443966297@qq.com}, Xue-ping Wang\footnote{Corresponding author. xpwang1@hotmail.com; fax: +86-28-84761393}\\
\emph{College of Mathematics and Software
Science, Sichuan Normal University,}\\
\emph{Chengdu, Sichuan 610066, People's Republic of China}}

\newcommand{\pp}[2]{\frac{\partial #1}{\partial #2}}
\date{}
\maketitle
\begin{quote}
{\bf Abstract}

This paper deals with the uniqueness of $L$-fuzzy sets in the representation of a given family of subsets of a nonempty set. It first shows a formula of the number of $L$-fuzzy sets whose collection of cuts coincides with a given family of subsets of a nonempty set, and then provides a necessary and sufficient condition under which such $L$-fuzzy set is unique. \\


\end{quote}
\section{Introduction}\label{intro}
In the classical paper \cite{Zadeh}, Zadeh introduced a notion of a fuzzy set of a set $X$ as a function from $X$ into [0,1]. In 1967, Goguen \cite{Goguen} gave a generalized version of the notion which is called an $L$-fuzzy set. Since then, $L$-fuzzy sets and structures have been widely studied. It is well-known that $L$-fuzzy mathematics attracts more and more interest in many
branches, for instance, algebraic theories including order-theoretic structures (see e.g., \cite{Gorjan14,Gorjan16,Saidi05,Seselja5}), automata and tree series
(see e.g., \cite{Borchar}) and theoretical computer science (see e.g., \cite{Chechik}). Among all the topics on $L$-fuzzy mathematics, the representation of a poset by an $L$-fuzzy set is very interesting, which in the case of a fuzzy set had been studied by \v{S}e\v{s}elja and Tepav\v{c}evi\'{c} \cite{Seselja2,Seselja98,Seselja4}, and by Jaballah and Saidi \cite{Jaballah} who investigated the characterization of all fuzzy sets of $X$ that can be identified with a given arbitrary family $C$ of subsets of $X$ together with a given arbitrary subset of [0,1], in particular, they gave necessary and sufficient conditions for both existence and uniqueness of such fuzzy sets. Later, Saidi and Jaballah investigated the problem of uniqueness under different considerations (see e.g., \cite{Saidi07,Saidi08,Saidi08F}). Also, Gorjanac-Ranitovi\'{c} and Tepav\v{c}evi\'{c} \cite{Gorjan} formulated a necessary and sufficient condition, under which for a given family of subsets $\mathcal{F}$ of a set $X$ and a fixed complete lattice $L$ there is an $L$-fuzzy set $\mu$ such that the collection of cuts of $\mu$ coincides with $\mathcal{F}$. Further, Jim\'{e}nez, Montes, \v{S}e\v{s}elja and Tepav\v{c}evi\'{c} \cite{Jorge} showed a necessary and sufficient condition, under which a collection of crisp up-sets (down-sets) of a poset $X$ consists of cuts of a lattice valued up-set (down-set). Therefore, a natural problem is: What is the condition under which the uniqueness of $L$-fuzzy set whose collection of cuts coincides with a given family of subsets of a nonempty set is guaranteed? Unfortunately, there are no results about the problem till now. In this paper, we shall discuss the condition under which such $L$-fuzzy set is unique.

The paper is organized as follows. For the sake of convenience, some
notions and previous results are given in Section 2. Section 3 first shows a formula of the number of $L$-fuzzy sets whose collection of cuts is equal to a given family of subsets of a nonempty set, and then provides a necessary and sufficient condition under which such $L$-fuzzy set is unique. Conclusions are drawn
in Section 4.
\section{Preliminaries}

We list some necessary notions and relevant properties from the classical order theory in the sequel. For more comprehensive presentation, see e.g., book \cite{Gratze}.

A poset is a structure $(P, \leq)$ where $P$ is a nonempty set and $\leq$ an ordering (reflexive, antisymmetric and transitive)
relation on $P$. A poset $(P, \leq_d)$ is a dual poset to the poset $(P, \leq)$, where $\leq_d$ is a dual ordering relation, defined by $x\leq_d y$ if and only if $y\leq x$. A sub-poset of a poset $(P, \leq)$ is a poset $(Q, \leq)$ where $Q$ is a nonempty subset of $P$ and $\leq$ on $Q$ is
restricted from $P$. A complete lattice is a poset $(L, \leq)$ in which every subset $M$ has the greatest lower bound, infimum,
meet, denoted by $\bigwedge M$, and the least upper bound, supremum, join, denoted by $\bigvee M$. A complete lattice $L$ possesses the top element
$1_L$ and the bottom element $0_L$.

In the following, we present some notations from the theory of $L$-fuzzy sets. More details about the relevant properties can be found e.g., in \cite{Gorjan, Jorge, Seselja3, Seselja4}.

An $L$-fuzzy set, Lattice-valued or $L$-valued set here is a mapping $\mu : X\rightarrow L$ from a nonempty set $X$ (domain) into a complete
lattice $(L, \wedge, \vee, 0_L, 1_L)$ (co-domain) (see \cite{Goguen}).

If $\mu : X\rightarrow L$ is an $L$-fuzzy set on $X$ then, for $p\in L$, the set \begin{equation*}\label{equ 00}\mu_p =\{x\in X\mid \mu(x)\geq p\}\end{equation*} is called the $p$-cut,
a cut set or simply a cut of $\mu$. Let $L^{\mu}$ be defined by \begin{equation}\label{equ 0}L^{\mu}=\{p\in L\mid p=\bigwedge B \mbox{ with } B\subseteq \mu(X)\},\end{equation}
where $\mu(X)=\{\mu(x)\mid x\in X\}$. Then $L^{\mu}$ is a complete lattice (see \cite{Jorge}).

We say that a mapping $f$: $L_1\rightarrow L_2$ from a complete lattice $L_1$ to a complete lattice $L_2$ preserves all infima if $f(\bigwedge_{s\in S}s)=\bigwedge_{s\in S}f(s)$ for all $S\subseteq L_1$, and it preserves the top element if $f(1_{L_{1}})=1_{L_{2}}$. Then the following two statements present some characterizations of the collection of cuts of an $L$-fuzzy set.

\begin{proposition}[\cite{Gorjan}]\label{pro}
Let $(L, \vee_L, \wedge_L)$ and $(L_{1}, \vee_{L_1}, \wedge_{L_1})$ be complete lattices and let $\phi: L\rightarrow L_1$ be the injection from $L$ to $L_1$ which maps the top element of $L$ to the top element of $L_1$, such that for all $x, y\in L, \phi(x\wedge_L y)= \phi(x)\wedge_{L_1} \phi(y)$. Let $\mu: X\rightarrow L$ be an $L$-fuzzy set. Let $L$-fuzzy set $\nu: X\rightarrow L_1$ be defined by $\nu(x)=\phi(\mu(x))$. Then the two $L$-fuzzy sets $\mu$ and $\nu$ have the same families of cuts and $\mu_p=\nu_{\phi(p)}$ for all $p\in L$.
\end{proposition}
\begin{theorem}[\cite{Gorjan}]\label{the 1}
Let $L$ be a fixed complete lattice. Necessary and sufficient conditions under which $\mathcal{F}\subseteq \mathcal{P}(X)$ is the collection of cut sets
of an $L$-fuzzy set with domain $X$ are: \\
\emph{(1)} $\mathcal{F}$ is closed under arbitrary intersections and contains $X$.\\
\emph{(2)} The dual poset of $\mathcal{F}$ under inclusion can be embedded into $L$, such that all infima and the top element are preserved under the embedding.
\end{theorem}

\section{Uniqueness of $L$-fuzzy sets}
In this section we shall investigate the condition for the uniqueness of $L$-fuzzy sets in the representation of
families of sets.

Let $(M, \leq)$ be a sub-poset of $(N, \leq)$. Define a mapping $\iota_{(M, N)}: M\rightarrow N$ by
\begin{equation}\label{e 1}\iota_{(M, N)} (x)=x\end{equation}
for all $x\in M$.
\begin{definition}\label{def 1}
\emph{Let $(L_1, \leq)$ be a sub-poset of a complete lattice $(L, \leq)$. The mapping $\iota_{(L_1, L)}$ is called an $\iota_{(L_1, L)}$-embedding if all infima and the top element are preserved under $\iota_{(L_1, L)}$.}
\end{definition}

Let $\mu_L=\{\mu_p\mid p\in L\}$ if $\mu: X\rightarrow L$ is an $L$-fuzzy set. Then we have the following theorem.
\begin{theorem}\label{theo 1}
Let $L$ be a complete lattice and let $\mu: X\rightarrow L$ be an $L$-fuzzy set. Then:\\
\emph{(a)} $(\mu_L, \supseteq) \cong (L^{\mu}, \leq)$;\\
\emph{(b)} $L^{\mu}$ can be embedded into $L$ by $\iota_{(L^{\mu}, L)}$-embedding;\\
\emph{(c)} $\mu=\iota_{(L^{\mu}, L)} \circ \nu$ where the mapping $\nu: X\rightarrow L^{\mu}$ is defined by $\nu(x)=\mu(x)$ for all $x\in X$, $\nu_{L^{\mu}}=\mu_L$ and $(L^{\mu})^{\nu}=L^{\mu}$.
\end{theorem}
\proof
(a) Let $\varphi: \mu_L\rightarrow L^{\mu}$ be defined by
\begin{equation}\label{equ 1}\varphi(\mu_p)=\bigwedge_{x\in \mu_p} \mu(x)\end{equation}
for all $\mu_p\in \mu_L$. It is easy to see that for all $p\in L$
\begin{equation}\label{000}\varphi(\mu_p)\geq p.\end{equation}

In what follows, we further prove that

\begin{equation}\label{001}\varphi(\mu_p)= p\end{equation}
for all $p\in L^{\mu}$.

Indeed, let $p\in L^{\mu}$. Then $p=\bigwedge_{B_p \subseteq \mu(X)} B_p$ by formula (\ref{equ 0}). If $B_p=\emptyset$ then $p=1_L$.
Clearly, $\varphi(\mu_{p})=\varphi(\mu_{1_L})=1_L=p$. Now, suppose $B_p \neq \emptyset$. Then from $B_p\subseteq \mu(X)=\{\mu(x)\mid x\in X\}$, we have $B_p\subseteq \{\mu(x)\in \mu(X)\mid \mu(x)\geq p\}$ since $b\in B_p$ implies $p\leq b$. Thus by formulas (\ref{equ 1}) and (\ref{000}),
\begin{equation*}\label{equ 2}p=\bigwedge_{B_p \subseteq \mu(X)} B_p \geq \bigwedge_{\mu(x)\geq p}
\mu(x)=\bigwedge_{x\in \mu_p}\mu(x)=\varphi(\mu_p)\geq p,\end{equation*} i.e., $\varphi(\mu_p)=p$. Therefore, $\varphi(\mu_p)= p$ for all $p\in L^{\mu}$. 

Formula (\ref{001}) means that $\varphi$ is surjective. Moreover, one can check that $\varphi$ is
injective by formulas (\ref{equ 1}) and (\ref{000}). Consequently, the mapping $\varphi$ is a bijection.

In what follows, we shall prove that both $\varphi$ and $\varphi^{-1}$ preserve the orders $\supseteq$ and $\leq$, respectively. In fact, if $\mu_p \subseteq \mu_q$ then obviously
$\varphi(\mu_p)=\bigwedge_{x\in \mu_p} \mu(x)\geq \bigwedge_{x\in \mu_q}\mu(x)=\varphi(\mu_q)$, i.e., $\varphi(\mu_p)\geq\varphi(\mu_q)$. At the same time, if $r, e\in L^{\mu}$ and $r\leq e$ then by the definition of a cut set, we have $\mu_r\supseteq \mu_e$. Thus, by formula (\ref{001}), $\varphi^{-1}(r)=\mu_r\supseteq \mu_e=\varphi^{-1}(e)$, i.e., $r\leq e$ implies $\varphi^{-1}(r)\supseteq \varphi^{-1}(e)$.

Therefore, $\varphi$ is an isomorphism from $(\mu_L, \supseteq)$ to $(L^{\mu}, \leq)$.

(b) From formula (\ref{equ 0}), $1_L=\bigwedge\emptyset \in L^{\mu}$ since $\emptyset\subseteq \mu(X)$. Thus by Definition \ref{def 1}, due to formulas (\ref{equ 0}) and (\ref{e 1}), $L^{\mu}$ can be embedded into $L$ by $\iota_{(L^{\mu}, L)}$-embedding.

(c) It is easy to see $\mu=\iota_{(L^{\mu}, L)} \circ \nu$. By (b), all the conditions
of Proposition \ref{pro} are fulfilled. Thus $\nu_{L^{\mu}}=\mu_L$. Moreover, by applying (a) to $L^{\mu}$-fuzzy set $\nu$, we know that $(\nu_{L^{\mu}}, \supseteq)\cong ((L^{\mu})^{\nu}, \leq)$. Then $(\mu_L, \supseteq)\cong ((L^{\mu})^{\nu}, \leq)$, which together with $(\mu_L, \supseteq)\cong (L^{\mu},\leq)$ yields
that $(L^{\mu},\leq)\cong ((L^{\mu})^{\nu}, \leq)$. Therefore, from (b), $(L^{\mu})^{\nu}=L^{\mu}$ since $(L^{\mu})^{\nu}$ is embedded into $L^{\mu}$ by $\iota_{((L^{\mu})^{\nu}, L^{\mu})}$-embedding.
\endproof
\begin{lemma}\label{lemm 1}
Let $\mathcal{F}$ be a family of some subsets of a nonempty set $X$ which is closed under intersections and contains $X$, and
let $L$ be a complete lattice. If $(\mathcal{F}, \supseteq)\cong (L_0, \leq)$ and $L_0$ can be embedded into $L$ by $\iota_{(L_0, L)}$-embedding, then:\\
\emph{(i)} there exists an $L_0$-fuzzy set $\nu: X\rightarrow L_0$ such that $\nu_{L_0}=\mathcal{F}$ and $L_{0}^{\nu}=L_0$;\\
\emph{(ii)} $\mu=\iota_{(L_0, L)}\circ \nu$ satisfies both $\mu_L=\mathcal{F}$ and $L^{\mu}=L_0$.
\end{lemma}
\proof By the hypotheses, $(\mathcal{F}, \supseteq)$ can be embedded into $L_0$, such that all infima and the top element are preserved under the embedding. Then by Theorem \ref{the 1}, there exists an $L_0$-fuzzy set $\nu: X\rightarrow L_0$ such that $\nu_{L_0}=\mathcal{F}$.
Thus by applying (a) of Theorem \ref{theo 1} to $\nu$, we have $(L_{0}^{\nu}, \leq) \cong (\nu_{L_0}, \supseteq) =(\mathcal{F}, \supseteq)\cong (L_0, \leq)$, i.e., $(L_{0}^{\nu}, \leq) \cong (L_0, \leq)$. Therefore,
$L_0 = L_{0}^{\nu}$ since $L_0 \supseteq L_{0}^{\nu}$ by formula (\ref{equ 0}).

Now, let $\mu=\iota_{(L_0, L)}\circ\nu$. Since $L_0$ can be embedded into $L$ by $\iota_{(L_0, L)}$-embedding, we know that all the conditions of
Proposition \ref{pro} are fulfilled. Then $\mu_L=\mathcal{F}$ since $\nu_{L_0}=\mathcal{F}$. Moreover, by formula (\ref{equ 0}) we have $L^{\mu}=L^{\nu}_{0}$ since $\mu(x)=\nu(x)$ for all $x\in X$ and $L_0$
can be embedded into $L$ by $\iota_{(L_0, L)}$-embedding. Therefore, $L^{\mu}=L_{0}$ since $L_{0}^{\nu}=L_0$.
\endproof

Let $\mathcal{F}$ be a family of some subsets of a nonempty set $X$ which is closed under intersections and contains $X$. Given a complete lattice $L$, we let $OI(L)$ denote the set of lattice isomorphisms from $L$ onto itself. Also, denote \[\begin{aligned}
&H_{(L_0, L_0, \mathcal{F})}=\{\mu \mid \mu: X\rightarrow L_0, L_{0}^{\mu}=L_0 \mbox { and } \mu_{L_0} =\mathcal{F}\}\mbox{ and }\\
&H_{(L, L_0, \mathcal{F})}=\{\mu \mid \mu: X\rightarrow L, L^{\mu}=L_0 \mbox { and } \mu_L =\mathcal{F}\}.\end{aligned}\]
Let
$$\mathcal{N}(L, \mathcal{F})=\{f\mid f: X\rightarrow L, f_{L}= \mathcal{F}\}$$
and $$\mathcal{S}(L, \mathcal{F})=\{ P\mid P\mbox{ can be embedded into }L\mbox{ by }\iota_{(P, L)}\mbox{-embedding and }(P, \leq)\cong (\mathcal{F}, \supseteq)\}.$$

Then, we have the following lemma.
\begin{lemma}\label{theo 3}
Let $\mathcal{F}$ be a family of some subsets of a nonempty set $X$ which is closed under intersections and contains $X$, and let $L$ be a complete lattice.
Then $$\mathcal{N}(L, \mathcal{F})=\bigcup_{L_0\in \mathcal{S}(L, \mathcal{F})}H_{(L, L_0, \mathcal{F})}.$$
\end{lemma}
\proof
Obviously, $\mathcal{N}(L, \mathcal{F})\supseteq\bigcup_{L_0\in \mathcal{S}(L, \mathcal{F})}H_{(L, L_0, \mathcal{F})}$. On the other hand, let $\mu \in \mathcal{N}(L, \mathcal{F})$. Then $\mathcal{F}=\mu_L$. Thus, by Theorem \ref{theo 1}, $(\mathcal{F}, \supseteq)=(\mu_L, \supseteq)\cong(L^{\mu}, \leq)$ and $L^{\mu}$ can be embedded into $L$ by $\iota_{(L^{\mu}, L)}$-embedding. Therefore $L^{\mu}\in \mathcal{S}(L, \mathcal{F})$ and $\mu \in H_{(L, L^{\mu}, \mathcal{F})}$, and we conclude $\mathcal{N}(L, \mathcal{F})\subseteq\bigcup_{L_0\in \mathcal{S}(L, \mathcal{F})}H_{(L, L_0, \mathcal{F})}$.
\endproof

\begin{theorem}\label{theo 2}
Under the assumptions of Lemma \ref{lemm 1}, we have:\\
\emph{(I)} $H_{(L, L_0, \mathcal{F})}=\{\beta\mid \beta=\iota_{(L_0, L)}\circ \mu, \mu\in H_{(L_0, L_0, \mathcal{F})}\}$;\\
\emph{(II)} Let $g\in H_{(L_0, L_0, \mathcal{F})}$. Then $f\in H_{(L_0, L_0, \mathcal{F})}$ if and only if $f=\eta\circ g$ for some $\eta\in OI(L_0)$.
\end{theorem}
\proof (I) We first note that $H_{(L_0, L_0, \mathcal{F})}\neq \emptyset$ by Lemma \ref{lemm 1}.
Then, by using Theorem \ref{theo 1} and Lemma \ref{lemm 1}, we have $H_{(L, L_0, \mathcal{F})}=\{\beta\mid \beta=\iota_{(L_0, L)}\circ \mu, \mu\in H_{(L_0, L_0, \mathcal{F})}\}$.

(II) Let $f\in H_{(L_0, L_0, \mathcal{F})}$. Define $\eta: L_0\rightarrow L_0$ by $$\eta(\bigwedge_{x\in B\subseteq X}g(x))=\bigwedge_{x\in B\subseteq X}f(x)$$
and $\eta_1: L_0\rightarrow L_0$ by $$\eta_1(\bigwedge_{x\in B\subseteq X}f(x))=\bigwedge_{x\in B\subseteq X}g(x).$$

The following proof is made in three parts:

 Part (1). Both $\eta$ and $\eta_1$ are mappings from $L_0$ to $L_0$.

Note that, from $f, g\in H_{(L_0, L_0, \mathcal{F})}$, we have $L_{0}^{f}=L_0$ and $L_{0}^{g}=L_0$, which together with (\ref{equ 2}) and (\ref{001}) imply that
\begin{equation}\label{e 2}t=\bigwedge_{x\in g_t}g(x)=\bigwedge_{x\in f_t}f(x)\end{equation} for any $t\in L_0$. Then
\begin{equation*}\label{003}\eta(t)=\eta(\bigwedge_{x\in g_t}g(x))=\bigwedge_{x\in g_t}f(x)\in L_0\mbox{ and }\eta_1(t)=\eta_1(\bigwedge_{x\in f_t}f(x))=\bigwedge_{x\in f_t}g(x)\in L_0\end{equation*} for all $t\in L_0$.

Therefore, in order to prove that $\eta$ is a mapping from $L_0$ to $L_0$, we just need to prove \begin{equation}\label{010}\eta(\bigwedge_{x\in B_1\subseteq X}g(x))=\eta(\bigwedge_{x\in B_2\subseteq X}g(x))\end{equation}
while $\bigwedge_{x\in B_1\subseteq X}g(x)= \bigwedge_{x\in B_2\subseteq X}g(x)$.

Now let $p=\bigwedge_{x\in B\subseteq X}g(x)$. Also, we can let $p=\bigwedge_{x\in g_p}g(x)$ by (\ref{e 2}). We shall show that \begin{equation}\label{e 3}\eta(\bigwedge_{x\in B\subseteq X}g(x))=\eta(\bigwedge_{x\in g_p}g(x))\end{equation} because formula (\ref{e 3}) implies (\ref{010}).

First, it is easy to see that $B\subseteq g_p$.
Then from $f, g\in H_{(L_0, L_0, \mathcal{F})}$, $f_{L_0}=g_{L_0}$. Thus $g_p\in g_{L_0}=f_{L_0}$, and again by formula (\ref{001}), there exists a unique element $w \in L^{f}_{0}$, i.e., $w \in L_0$ such that \begin{equation}\label{004}f_w =g_p\end{equation} since $L_0=L^{f}_{0}$. So that, by (\ref{e 2}), \begin{equation}\label{008}\eta(\bigwedge_{x\in g_p}g(x))=\bigwedge_{x\in g_p}f(x)=\bigwedge_{x\in f_w}f(x)=w.\end{equation}

On the other hand, by the definition of $\eta$ and $L_0=L^{f}_0$, we have $\eta(\bigwedge_{x\in B\subseteq X}g(x))=\bigwedge_{x\in B\subseteq X}f(x)\in L_0$.
Let $\bigwedge_{x\in B\subseteq X}f(x)=r$, i.e., \begin{equation}\label{009}\eta(\bigwedge_{x\in B\subseteq X}g(x))=r.\end{equation} Then \begin{equation}\label{005}B\subseteq f_r,\end{equation} and $r\geq w$ since $B\subseteq g_p$. Thus \begin{equation}\label{006}f_r\subseteq f_w.\end{equation} Using (\ref{e 2}), we have $\bigwedge_{x\in B\subseteq X}f(x)=r=\bigwedge_{x\in f_r}f(x)$. Clearly $f_r \in f_{L_0}=g_{L_0}$. Again, by formula (\ref{001}), there exists a unique element $q \in L^{g}_{0}$, i.e., $q \in L_0$ such that \begin{equation}\label{007}f_r =g_q\end{equation} since $L_0=L^{g}_{0}$. Thus, by formulas (\ref{005}), (\ref{006})and (\ref{004}), $B\subseteq g_q \subseteq g_p$.
Hence, by (\ref{e 2}) we know that  $$p=\bigwedge_{x\in B\subseteq X}g(x)\geq \bigwedge_{x\in g_q}g(x)=q \geq \bigwedge_{x\in g_p}g(x)=p$$ i.e., $p=q$.
Then $g_p =g_q$, and which together with (\ref{004}) and (\ref{007}) means that $f_r=f_w$. Thus, by (\ref{001}),
we conclude that $r=w$, i.e., \begin{equation*}\eta(\bigwedge_{x\in B\subseteq X}g(x))=\eta(\bigwedge_{x\in g_p}g(x))\end{equation*}
by (\ref{008}) and (\ref{009}). This completes the proof of formula (\ref{e 3}).

Consequently $\eta$ is a mapping from $L_0$ to itself.

Similarly, $\eta_1$ is also a mapping from $L_0$ to itself.

Part (2). Both $\eta$ and $\eta_1$ preserve the order.

Assume that $t_1, t_2 \in L_0$ and $t_1 \leq t_2$. Then $g_{t_1}\supseteq g_{t_2}$ by the definition of a cut set. Thus by (\ref{e 2}), $\eta(t_1)=\bigwedge_{x\in g_{t_1}} f(x) \leq \bigwedge_{x\in g_{t_2}} f(x)=\eta(t_2)$,
and which means that $\eta$ preserves the order. Also, we can similarly check that $\eta_1$ preserves the order.

Part (3). The mapping $\eta$ is a bijection and $\eta_1=\eta^{-1}$.

Let $k\in L_0$. Then from (\ref{e 2}),
$$\eta_1\circ \eta(k)=\eta_1\circ \eta(\bigwedge_{x\in g_{k}} g(x))=\eta_1(\bigwedge_{x\in g_{k}} f(x))=\bigwedge_{x\in g_{k}} g(x)=k.$$
Thus $\eta_1\circ \eta$ is an identity mapping on $L_0$. Similarly, we can check that $\eta\circ\eta_1$ is also an identity mapping on $L_0$.
Therefore, the mapping $\eta$ is a bijection and $\eta_1=\eta^{-1}$.

From Parts (1), (2) and (3), we conclude that $\eta\in OI(L_0)$ and $f=\eta\circ g$.

Conversely, let $\eta\in OI(L_0)$ and $\eta\circ g=f$. We shall show that $f\in H_{(L_0, L_0, \mathcal{F})}$.

Since $f$ is an $L_0$-fuzzy set on $X$, we need to prove that $f_{L_0}=\mathcal{F}$ and $L^f_0=L_0$.

First, we prove $f_{L_0}=\mathcal{F}$, i.e., prove
that $f_{L_0}=g_{L_0}$ since $g\in H_{(L_0, L_0, \mathcal{F})}$.

From $\eta \in OI(L_0)$,
\[\begin{aligned}
f_{\eta(s)}&=\{x\in X \mid \eta\circ g(x)\geq \eta(s)\}\\
           &= \{x\in X \mid g(x)\geq \eta^{-1}(\eta(s))=s\}\\
           &=g_s
\end{aligned}\] for all $s\in L_0$. Therefore $g_{L_0}\subseteq f_{L_0}$.
Similarly, we can prove that $g_{L_0}\supseteq f_{L_0}$ since $g=\eta^{-1}\circ f$ and $\eta^{-1}\in OI(L_0)$. Hence, \begin{equation}\label{011}f_{L_0}=g_{L_0}.\end{equation}

Secondly, we shall prove $L^f_0=L_0$.

By Theorem \ref{theo 1}, we know that $(f_{L_0}, \supseteq)\cong (L_{0}^{f}, \leq)$ and $(g_{L_0}, \supseteq)\cong (L_{0}^{g}, \leq)$. Thus, by (\ref{011}), $(L_{0}^{f}, \leq)\cong(L_{0}^{g}, \leq)$. Note that $g\in H_{(L_0, L_0, \mathcal{F})}$ implies that $(L_{0}^{g}, \leq)=(L_0, \leq)$. Therefore $(L_{0}^{f}, \leq)\cong(L_{0}, \leq)$.
This follows that $L_{0}^{f}=L_0$ since $L_{0}^{f}\subseteq L_0$.

Finally, $f\in H_{(L_0, L_0, \mathcal{F})}$.
\endproof
From Theorem \ref{theo 2}, we easily deduce the following consequence.

\begin{corollary}\label{cor1}
Under the assumptions of Lemma \ref{lemm 1}. Let $g\in H_{(L_0, L_0, \mathcal{F})}$. Then
\[\begin{aligned}
&H_{(L_0, L_0, \mathcal{F})}=\{\eta\circ g\mid \eta \in OI(L_0)\} \mbox{ and }\\
&H_{(L, L_0, \mathcal{F})}=\{\iota_{(L_0, L)}\circ\eta\circ g\mid \eta \in OI(L_0)\}.
\end{aligned}\]
\end{corollary}

Given any set $A$, we denote its cardinality by $|A|$. Then from Corollary \ref{cor1} and Lemma \ref{theo 3}, we have:
\begin{theorem}\label{cor2}
Let $\mathcal{F}$ be a family of some subsets of a nonempty set $X$ which is closed under intersection and contains $X$, and let $L$ be a complete lattice.
Then $$|\mathcal{N}(L, \mathcal{F})|=|\mathcal{S}(L, \mathcal{F})| |OI(\mathcal{F})|.$$
\end{theorem}
\proof
From Lemma \ref{theo 3}, this result obviously holds if $\mathcal{S}(L, \mathcal{F})=\emptyset$. Let $L_0 \in \mathcal{S}(L, \mathcal{F})$. Then $(L_0, \leq)\cong (\mathcal{F}, \supseteq)$, which means that \begin{equation}\label{e 6}|OI(\mathcal{F})|=|OI(L_0)|.\end{equation} Now, by Corollary \ref{cor1}, we have \begin{equation}\label{e 4}|H_{(L, L_0, \mathcal{F})}|=|H_{(L_0, L_0, \mathcal{F})}|.\end{equation}

Let $g\in H_{(L_0, L_0, \mathcal{F})}$. Suppose $f_1, f_2 \in H_{(L_0, L_0, \mathcal{F})}$. Then by Corollary \ref{cor1}, there exist two mapping $\eta_1, \eta_2 \in OI(L_0)$ such that $f_1=\eta_1\circ g$
and $f_2=\eta_2\circ g$. Thus $f_1\neq f_2$ if and only if $\eta_1 \neq\eta_2$, which results in \begin{equation}\label{e 5}|H_{(L_0, L_0, \mathcal{F})}|=|OI(L_0)|.\end{equation}

Let $L_1\in \mathcal{S}(L, \mathcal{F})$ with $L_1\neq L_0$. Then it is clear that $H_{(L, L_1, \mathcal{F})}\bigcap H_{(L, L_0, \mathcal{F})}=\emptyset$, and $|H_{(L, L_1, \mathcal{F})}|=|H_{(L, L_0, \mathcal{F})}|$ since $(L_1, \leq)\cong (L_0, \leq)$.
Therefore, from Lemma \ref{theo 3} and formulas (\ref{e 6}), (\ref{e 4}) and (\ref{e 5}), we have $|\mathcal{N}(L, \mathcal{F})|=|\mathcal{S}(L, \mathcal{F})| |OI(\mathcal{F})|$.
\endproof

The following example illustrates Theorem \ref{cor2}.
\begin{example}\label{exa 1}
\emph{Let us consider a set $X=\{a, b, c\}$ and the complete lattice $(L, \leq)$ represented in Fig.1. Let
$\mathcal{F}=\{\emptyset, \{b\}, \{a, b\}, \{b, c\}, \{a, b, c\}\}$ be a family of subsets of $X$.}
\par\noindent\vskip50pt
 \begin{minipage}{11pc}
\setlength{\unitlength}{0.75pt}\begin{picture}(600,100)
\put(270,40){\circle{4}}\put(275,30){\makebox(0,0)[l]{\footnotesize $0_L$}}
\put(235,75){\circle{4}}\put(230,65){\makebox(0,0)[l]{\footnotesize $q$}}
\put(270,75){\circle{4}}\put(275,65){\makebox(0,0)[l]{\footnotesize $r$}}
\put(305,75){\circle{4}}\put(310,65){\makebox(0,0)[l]{\footnotesize $p$}}
\put(235,110){\circle{4}}\put(230,120){\makebox(0,0)[l]{\footnotesize $s$}}
\put(305,110){\circle{4}}\put(310,120){\makebox(0,0)[l]{\footnotesize $t$}}
\put(270,145){\circle{4}}\put(275,150){\makebox(0,0)[l]{\footnotesize $1_L$}}
\put(235,77){\line(0,1){31}}
\put(268.5,41.5){\line(-1,1){32}}
\put(271.5,41.5){\line(1,1){32}}
\put(270,42){\line(0,1){31}}
\put(271.5,76.5){\line(1,1){32}}
\put(268.5,76.5){\line(-1,1){32}}
\put(305,77){\line(0,1){31}}
\put(236.5,111.5){\line(1,1){32}}
\put(303.5,111.5){\line(-1,1){32}}
\put(270,5){$L$}

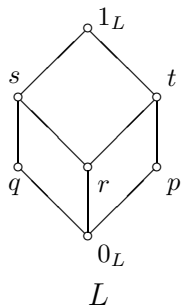
\put(175,-15){\emph{Fig.1 Hasse diagram of $L$}}
  \end{picture}
  \end{minipage}
$$ $$

\emph{It is clear that $\mathcal{S}(L, \mathcal{F})=\{(\{0_L, p, r, t, 1_L\}, \leq), (\{0_L, q, r, s, 1_L\}, \leq)\}$ and $|OI(\mathcal{F})|=2$.
Therefore, from Theorem \ref{cor2}, $|\mathcal{N}(L, \mathcal{F})|=4$. On the other hand, all the $L$-fuzzy sets of $\mathcal{N}(L, \mathcal{F})$ are}
\begin{table}[htbp]
\centering
\begin{tabular}{|c|c|c|c|}
  \hline
 $X$ & $a$ & $b$ & $c$ \\
 \hline
$\delta$ & $r$ & $t$ & $p$ \\
  \hline
\end{tabular} \ \ \ \ \ \ \ \ \ \ \begin{tabular}{|c|c|c|c|}
  \hline
  $X$ & $a$ & $b$ & $c$ \\
   \hline
 $\gamma$ & $p$ & $t$ & $r$ \\
  \hline
\end{tabular}
$$
\begin{tabular}{|c|c|c|c|}
  \hline
$X$ & $a$ & $b$ & $c$ \\
   \hline
 $\beta$ & $q$ & $s$ & $r$ \\
  \hline
\end{tabular} \ \ \ \ \ \ \ \ \ \ \ \begin{tabular}{|c|c|c|c|}
  \hline
$X$ & $a$ & $b$ & $c$ \\
  \hline
$\alpha$ & $r$ & $s$ & $q$ \\
  \hline
\end{tabular}
$$
 \end{table}\\
\emph{where we consider a tabular representation for the $L$-fuzzy set on $X$, respectively.}
\end{example}

The following theorem follows immediately from Theorem \ref{cor2}.
\begin{theorem}\label{theo 4}
Let $\mathcal{F}$ be a family of some subsets of a nonempty set $X$ which is closed under intersection and contains $X$, and let $L$ be a complete lattice.
There exists a unique $L$-fuzzy $\mu$ on $X$ such that $\mathcal{F}=\mu_{L}$ if and
only if $|\mathcal{S}(L, \mathcal{F})|=|OI(\mathcal{F})|=1$.
\end{theorem}

\section{Conclusions}

This contribution gave a necessary and sufficient condition under which $L$-fuzzy sets whose collection of cuts equals a given family of subsets of a nonempty set are unique. Using Theorem 12 in \cite{Jorge}, one can verify that all our results made for $L$-fuzzy sets could be applied to $L$-fuzzy up-sets ($L$-fuzzy down-sets) since an $L$-fuzzy up-set ($L$-fuzzy down-set) is just a particularity of the concept of an $L$-fuzzy set.

\section*{Acknowledgments}
The authors thank the referees for their valuable comments and suggestions.

\end{document}